\documentclass[oneside,english]{amsart}
\usepackage[latin1]{inputenc}
\usepackage{geometry}
\usepackage{amssymb}

\makeatletter
 \theoremstyle{plain}    
 \newtheorem{thm}{Theorem}[section]
 \numberwithin{equation}{section} 
 \numberwithin{figure}{section} 
 \theoremstyle{plain}
 \theoremstyle{definition}
 \newtheorem{defn}[thm]{Definition}
 \theoremstyle{plain}    
 \newtheorem{lem}[thm]{Lemma} 
 \theoremstyle{plain}    
 \newtheorem{prop}[thm]{Proposition} 
 \theoremstyle{plain}    
 \newtheorem{cor}[thm]{Corollary} 
 \theoremstyle{definition}
  \newtheorem{example}[thm]{Example}

\newcommand{\Graph}{\operatorname{Graph}}

\newcommand{\MWGs}{ \mathcal{G}_{i}=(G,(K_{v}^{i})_{v\in V},(\phi_{e}^{i})_{e\in E}) }

\usepackage{babel}
\makeatother
\begin{document}

\title{Mauldin-Williams graphs, Morita Equivalence and Isomorphisms}

\author{Marius Ionescu}

\begin{abstract}
We describe a method for associating some non-self-adjoint algebras
to Mauldin-Williams graphs and we study the Morita equivalence and
isomorphism of these algebras.

We also investigate the relationship between the Morita equivalence
and isomorphism class of the $C^{\ast}$-correspondences associated
with Mauldin-Williams graphs and the dynamical properties of the Mauldin-Williams
graphs.
\end{abstract}
\maketitle

\section{Introduction}

In this note we follow the notation from \cite{John}. By a \emph{Mauldin-Williams}
\emph{graph} (see \cite{Mauldin-Williams}), we mean a system $\mathcal{G}=(G,\{ T_{v},\rho_{v}\}_{v\in V},\{\phi_{e}\}_{e\in E})$,
where $G=(V,E,r,s)$ is a graph with a finite set of vertices $V$,
a finite set of edges $E$, a \emph{range} map $r$ and a \emph{source}
map $s$, and where $\{ T_{v},\rho_{v}\}_{v\in V}$ and $\{\phi_{e}\}_{e\in E}$
are families such that:

\begin{enumerate}
\item Each $T_{v}$ is a compact metric space with a prescribed metric $\rho_{v}$,
$v\in V$.
\item For $e\in E$, $\phi_{e}$ is a continuous map from $T_{r(e)}$ to
$T_{s(e)}$ such that\[
c_{1}\rho_{r(e)}(x,y)\le\rho_{s(e)}(\phi_{e}(x),\phi_{e}(y))\leq c\rho_{r(e)}(x,y)\]
 for some constants $c_{1},c$ satisfying $0<c_{1}\le c<1$ (independent
of $e$) and all $x,y\in T_{r(e)}$. 
\end{enumerate}
We shall assume, too, that the source map $s$ and the range map $r$
are surjective. Thus, we assume that there are no sinks and no sources
in the graph $G$.
\newcommand{\twolines}[2]{{\genfrac{}{}{0pt}{}{#1}{#2}}}

In \cite{John} we associated to a Mauldin-Williams graph $\mathcal{G}=(G,\{ T_{v},\rho_{v}\}_{v\in V},\{\phi_{e}\}_{e\in E})$
a so-called $C^{*}$-\emph{correspondence} $\mathcal{X}$ over the
$C^{*}$-algebra $A=C(T)$, where $T=\coprod_{v\in V}T_{v}$ is the
disjoint union of the spaces $T_{v},v\in V$, as follows. Let $E\times_{G}T=\{(e,x)\:|\: x\in T_{r(e)}\}$.
Then, by our finiteness assumptions, $E\times_{G}T$ is a compact
space. We set $\mathcal{X}=C(E\times_{G}T)$ and view $\mathcal{X}$
as a $C^{\ast}$-correspondence over $C(T)$ via the formulae:\begin{eqnarray*}
\xi\cdot a(e,x) & := & \xi(e,x)a(x),\\
a\cdot\xi(e,x) & := & a\circ\phi_{e}(x)\xi(e,x)\end{eqnarray*}
and\[
\langle\xi,\eta\rangle_{A}(x):=\sum_{{\genfrac{}{}{0pt}{}{e\in E}{x\in T_{r(e)}}}}\overline{\xi(e,x)}\eta(e,x),\]
 where $a\in C(T)$ and $\xi,\eta\in C(E\times_{G}T)$. With these
data we can form the tensor algebra $\mathcal{T}_{+}(\mathcal{X})$
as prescribed in \cite{Muhly-Baruc-tensoalgebras} and \cite{Muhly-Baruc-Moritaequiv}.
Our main result is:

\begin{thm}
\label{thm:Main}For $i=1,2$, let $\mathcal{G}_{i}=(G_{i},(K_{v}^{i})_{v\in V_{i}},(\phi_{e}^{i})_{e\in E_{i}})$
be two Mauldin-Williams graphs. Let $A_{i}=C(K^{i})$ and let $\mathcal{X}_{i}$
be the associated $C^{*}$-algebras and $C^{*}$-correspondences.
Then the following are equivalent:
\begin{enumerate}
\item \label{enu:first}$\mathcal{T}_{+}(\mathcal{X}_{1})$ is strongly
Morita equivalent to $\mathcal{T}_{+}(\mathcal{X}_{2})$ in the sense
of \cite{Blecher-Muhly-Paulsen}.
\item \label{enu:second}$\mathcal{X}_{1}$ and $\mathcal{X}_{2}$ are strongly
Morita equivalent in the sense of \cite{Muhly-Baruc-Moritaequiv}. 
\item \label{enu:third}$\mathcal{X}_{1}$ and $\mathcal{X}_{2}$ are isomorphic
as $C^{\ast}$-correspondeces.
\item \label{enu:fourth}$\mathcal{T}_{+}(\mathcal{X}_{1})$ is completely
isometrically isomorphic to $\mathcal{T}_{+}(\mathcal{X}_{2})$.
\end{enumerate}
\end{thm}
We find this result especially remarkable in light of Theorem 2.3
from \cite[Theorem 1.1]{John} (see also Section 4.2 from \cite{Pinzari-Watatami-Y}),
which states that the Cuntz-Pimsner algebra, $\mathcal{O}(\mathcal{X})$,
which is the $C^{\ast}$-envelope of the tensor algebra $\mathcal{T}_{+}(\mathcal{X})$,
depends only of the structure of the underlying graph.  In particular,
our results lead to examples of different non-self-adjoint algebras
which are not completely isometrically isomorphic, but have the same
$C^{*}$-envelope, namely $\mathcal{O}_{n}$.

To understand further the relationship between the tensor algebra
and the Mauldin-Williams graph, we study the isomorphism class of
our $C^{*}$-correspondences and tensor algebras in terms of the dynamics
of the Mauldin-Williams graph. Roughly, we find that two $C^{*}$-correspondences
associated to two Mauldin-Williams graphs, $(G_{i},(K_{v}^{(i)})_{v\in V_{i}},(\phi_{e}^{i})_{e\in E_{i}})$,
$i=1,2$ are isomorphic if the maps $\phi_{e}^{1}$ and $\phi_{e}^{2}$
are locally conjugate in a sense that will be made precise later.

\section{Non-self-adjoint Algebras Associated with Mauldin-Williams Graphs}

\begin{defn}
An \emph{invariant list} associated with a Mauldin-Williams graph
$\mathcal{G}=(G,\{ T_{v},\rho_{v}\}_{v\in V},\{\phi_{e}\}_{e\in E}\}$
is a family $(K_{v})_{v\in V}$ of compact sets, such that $K_{v}\subset T_{v}$,
for all $v\in V$ and such that\begin{align*}
K_{v} & =\bigcup_{e\in E,s(e)=v}\phi_{e}(K_{r(e)}).\end{align*}
Since each $\phi_{e}$ is a proper contraction, $\mathcal{G}$ has
a unique invariant list (see \cite[Theorem 1]{Mauldin-Williams}).
We set $T:=\bigcup_{v\in V}T_{v}$ and $K:=\bigcup_{v\in V}K_{v}$
and we call $K$ the \emph{invariant set} of the Mauldin-Williams
graph.
\end{defn}
In the particular case when we have one vertex $v$ and $n$ edges,
i.e. in the setting of an \emph{iterated function system,} the invariant
set is the unique compact subset $K:=K_{v}$ of $T=T_{v}$ such that\[
K=\phi_{1}(K)\cup\cdots\cup\phi_{n}(K).\]
Note that the $\ast$-homomorphism $\Phi:A\rightarrow\mathcal{L}(\mathcal{X})$,
$(\Phi(a)\xi)(e,x)=a\circ\phi_{e}(x)\xi(e,x)$, which gives the left
action of the $C^{\ast}$-correspondence associated to a Mauldin-Williams
graph, is faithful if and only if $K=T$. In this note we assume that
$T$ equals the invariant set $K$.

Kajiwara and Watatani have proved in \cite[Lemma 2.3]{KWKMSonCP}
that, if the contractions are proper, the invariant set of an iterated
function system has no isolated point. Their proof can be easily generalized
to the invariant set of a Mauldin-Williams graph. Hence $K$ has no
isolated points.

For a $C^{*}$-correspondence $\mathcal{X}$ over a $C^{*}$-algebra
$A$, the (full) \emph{Fock space} over $\mathcal{X}$ is\[
\mathcal{F}(\mathcal{X})=A\oplus\mathcal{X}\oplus\mathcal{X}^{\otimes2}\oplus\cdots.\]
 We write $\Phi_{\infty}$ for the left action of $A$ on $\mathcal{F}(\mathcal{X})$,
$\Phi_{\infty}(a)=\mbox{diag}(a,\Phi^{(1)}(a),\Phi^{(2)}(a),\cdots)$,
where $\Phi^{(n)}$ is the left action of $A$ on $\mathcal{X}^{\otimes n}$
($\Phi^{(1)}=\Phi$, the left action of $A$ on $\mathcal{X}$). For
$\xi\in\mathcal{X}$, the creation operator determined by $\xi$ is
defined by the formula $T_{\xi}(\eta)=\xi\otimes\eta$, for all $\eta\in\mathcal{F}(\mathcal{X})$.

\begin{defn}
The \emph{tensor algebra} of $\mathcal{X}$, denoted by $\mathcal{T}_{+}(\mathcal{X})$,
is the norm closed subalgebra of $\mathcal{L}(\mathcal{F}(\mathcal{X}))$
generated by $\Phi_{\infty}(A)$ and the creation operators $T_{\xi}$,
for $\xi\in\mathcal{X}$ (see \cite{Muhly-Baruc-tensoalgebras} and
\cite{Muhly-Baruc-Moritaequiv}). The $C^{*}$-algebra generated by
$\mathcal{T}_{+}(\mathcal{X})$ is denoted by $\mathcal{T}(\mathcal{X})$
and it is called the \emph{Toeplitz algebra} of the $C^{\ast}$-correspondence
$\mathcal{X}$.

We may regard each finite sum $\sum_{n=0}^{N}\mathcal{X}^{\otimes n}$
as a subspace of $\mathcal{F}(\mathcal{X})$ and we may regard $\mathcal{L}(\sum_{n=0}^{N}\mathcal{X}^{\otimes n})$
as a subalgebra of $\mathcal{L}(\mathcal{F}(\mathcal{X}))$ in the
obvious way. Let $B$ be the $C^{*}$-subalgebra of $\mathcal{L}(\mathcal{F}(\mathcal{X}))$
generated by all the $\mathcal{L}(\sum_{n=0}^{N}\mathcal{X}^{\otimes n})$
as $N$ ranges over the non-negative integers. Then $\mathcal{T}(\mathcal{X})\subset M(B)$,
the multiplier algebra of $B$. The Cuntz-Pimsner algebra $\mathcal{O}(\mathcal{X})$
is defined to be the image of $\mathcal{T}(\mathcal{X})$ in the corona
algebra $M(B)/B$ (see \cite{Muhly-Baruc-tensoalgebras} and \cite{Pimnser}). 

By a \emph{homomorphism} from an $A_{1}-B_{1}$ $C^{\ast}$-correspondence
$\mathcal{X}_{1}$, \emph{to an $A_{2}-B_{2}$ $C^{\ast}$-correspondence
$\mathcal{X}_{2}$} we mean a triple $(\alpha,V,\beta)$, where $\alpha:A_{1}\rightarrow A_{2}$,
$\beta:B_{1}\rightarrow B_{2}$ are $C^{\ast}$-homomorphisms and
$V:\mathcal{X}_{1}\rightarrow\mathcal{X}_{2}$ is a linear map such
that $V(a\xi b)=\alpha(a)V(\xi)\beta(b)$ and such that $\langle V(\xi),V(\eta)\rangle_{B_{2}}=\beta(\langle\xi,\eta\rangle_{B_{1}})$
(see \cite[Section 1]{Muhly-Baruc-Moritaequiv}). When $A_{1}=A_{2}$
and $B_{1}=B_{2}$, we will consider 
\newcommand{\Aut}{\operatorname{Aut}}
$\alpha\in\Aut(A_{1})$ and $\beta\in\Aut(B_{1})$. This, then, forces
$V$ to be isometric. If $V$ is also surjective, we shall say that
$V$ is a \emph{correspondence isomorphism} \emph{over} $(\alpha,\beta)$.
If, moreover, $A_{1}=B_{1}$ and $\alpha=\beta$, we say that $V$
is a \emph{correspondence isomorphism over $\alpha$.}
\end{defn}
A central concept for our work in this note is the \emph{strong Morita
equivalence} for $C^{*}$-correspondences defined in \cite[Definition 2.1]{Muhly-Baruc-Moritaequiv},
which we review here.
\newcommand{\mequiv}[3]{#1\overset{\operatorname{SME}}{\sim}_{#3}#2}

\begin{defn}
If $\mathcal{X}$ is a $C^{*}$-correspondence over a $C^{\ast}$-algebra
$A$, and $\mathcal{Y}$ is a $C^{*}$-correspondence over a $C^{\ast}$-algebra
$B$, we say that $\mathcal{X}$ and $\mathcal{Y}$ are \emph{strongly
Morita equivalent} if $A$ and $B$ are strongly Morita equivalent
via an $A$-$B$ equivalence bimodule $\mathcal{Z}$ (in which case
we write $\mequiv{A}{B}{\mathcal{Z}}$), for which there is an $A$-$B$
correspondence isomorphism $(id,W,id)$ from $\mathcal{Z}\otimes_{B}\mathcal{Y}$
onto $\mathcal{X}\otimes_{A}\mathcal{Z}$. This means, in particular,
that $W(a\xi b)=aW(\xi)b$ for all $a\in A,b\in B$ and $\xi\in\mathcal{Z}\otimes_{B}\mathcal{Y}$
and that $\langle W(\xi),W(\eta)\rangle_{B}=\langle\xi,\eta\rangle_{B}$.
\end{defn}
We say that a $C^{*}$-correspondence $\mathcal{X}$ over a $C^{*}$-algebra
$A$ is \emph{aperiodic} if: for all $n\geq1,\;\mbox{for all}\;\xi\in\mathcal{X}^{\otimes n}$
and for all hereditary subalgebras $B\subseteq A$, we have\[
\inf\left\{ \left\Vert \Phi^{(n)}(a)\xi a\right\Vert \;|\; a\geq0,a\in B,\Vert a\Vert=1\right\} =0.\]
It was proved in \cite[Theorem 3.2, Theorem 3.5]{Muhly-Baruc-Moritaequiv}
that if $\mathcal{X}$ and $\mathcal{Y}$ are strongly Morita equivalent,
then $\mathcal{T}_{+}(\mathcal{X})$ and $\mathcal{T}_{+}(\mathcal{Y})$
(respectively $\mathcal{T}(\mathcal{X})$ and $\mathcal{T}(\mathcal{Y})$
, $\mathcal{O}(\mathcal{X})$ and $\mathcal{O}(\mathcal{Y})$) are
strongly Morita equivalent. Also, if $\mathcal{X}$ and $\mathcal{Y}$
are aperiodic $C^{*}$-correspondences over the $C^{*}$-algebras
$A$ and $B$, respectively, and if $\mathcal{T}_{+}(\mathcal{X})$
and $\mathcal{T}_{+}(\mathcal{Y})$ are strongly Morita equivalent
in the sense of \cite{Blecher-Muhly-Paulsen}, then $\mathcal{X}$
and $\mathcal{Y}$ are strongly Morita equivalent (see \cite[Theorem 7.2]{Muhly-Baruc-Moritaequiv}).

To study the aperiodicity and strong Morita equivalence of $C^{\ast}$-correspondences
associated to Mauldin-Williams graphs, we need the following lemma
which gives an equivalent description of when a $C^{\ast}$-correspondence
is aperiodic.

\begin{lem}
\label{lem:aperiodic}(\cite[Lemma 5.2]{Muhly-Baruc-Moritaequiv}).
The $C^{*}$-correspondence $\mathcal{X}$ is aperiodic if and only
if given $a_{0}\in A$, $a_{0}\geq0$, $\xi^{k}\in\mathcal{X}^{\otimes k}$,
$1\leq k\leq n$ and $\varepsilon>0,$ there is an $x$ in the hereditary
subalgebra $\overline{a_{0}Aa_{0}}$, with $x\geq0$ and $\Vert x\Vert=1$,
such that\[
\Vert xa_{0}x\Vert>\Vert a_{0}\Vert-\varepsilon\]
 and

\[
\Vert\Phi^{(k)}(x)\xi^{k}x\Vert<\varepsilon\;\mbox{for}\;1\leq k\leq n.\]

\end{lem}
For a directed graph $G=(V,E,r,s)$ and for $k\geq2$, we define\[
E^{k}:=\{\alpha=(\alpha_{1},\cdots,\alpha_{k})\;:\;\alpha_{i}\in E\;\mbox{and}\; r(\alpha_{i})=s(\alpha_{i+1}),i=1,\cdots,k-1\}\]
 to be the set of paths of length $k$ in the graph $G$. We define
also the infinite path space to be\[
E^{\infty}:=\{(\alpha_{i})_{i\in\mathbb{N}}\;:\;\alpha_{i}\in E\;\mbox{and}\; r(\alpha_{i})=s(\alpha_{i+1})\;\mbox{for all}\: i\in\mathbb{N}\}\]
 For $\alpha\in E^{k}$, we write $\phi_{\alpha}=\phi_{\alpha_{1}}\circ\cdots\circ\phi_{\alpha_{k}}$.

\begin{prop}
\label{pro:aperiodicity}Let $\mathcal{G}=(G,(K_{v})_{v\in V},(\phi_{e})_{e\in E})$
be a Mauldin-Williams graph with the invariant set $K$. Let $A=C(K)$
be the associated $C^{*}$-algebra and let $\mathcal{X}$ be the associated
$C^{*}$-correspondence. Then the $C^{*}$-correspondence $\mathcal{X}$
is aperiodic. 
\end{prop}
\begin{proof}
Note that $\phi_{\alpha}:K_{r(\alpha)}\rightarrow K_{s(\alpha)}$,
with $\alpha\in E^{k}$ and $k\in\mathbb{N}$, has a fixed point if
and only if $r(\alpha)=s(\alpha)$, i.e. $\alpha$ is a cycle in the
graph $G$.

Fix $n_{0}\in\mathbb{N}$, choose $k\in\mathbb{N},1\le k\le n_{0}$;
let $a_{0}\in A$ with $a_{0}\geq0$; let $\xi^{k}\in\mathcal{X}^{\otimes k}$
and let $\varepsilon>0$. We verify the criterion in Lemma \ref{lem:aperiodic}
first when $n_{0}=k=1$.

Without loss of generality, we assume that $\left\Vert a_{0}\right\Vert =1$.
Then we can find $t_{0}\in K$ such that $|a_{0}(t_{0})|\geq1-\varepsilon$
and $t_{0}$ is not a fixed point for any $\phi_{e}$ , $e\in E$.
Let $v_{0}\in V$ be such that $t_{0}\in K_{v_{0}}$. Choose $\delta_{1}>0$
such that $B(t_{0},\delta_{1})\subset K_{v_{0}}$ and $B(\phi_{e}(t_{0}),\delta_{1})\cap B(t_{0},\delta_{1})=\emptyset$
for all $e\in E$ for which $r(e)=v_{0}$. Let\[
\delta_{2}:=\begin{cases}
\min\{\rho_{v_{0}}(t_{0},t)\;|\; a_{0}(t)=0\}, & \mbox{if}\;\{ t\in K_{v_{0}}\;:\; a_{0}(t)=0\}\ne\emptyset\\
\delta_{1}, & \mbox{otherwise}.\end{cases}\]
Set $\delta=\min\{\delta_{1},\delta_{2}\}$ and let $x\in A$, $x\geq0$
be such that

\[
x(t)=\left\{ \begin{array}{ccl}
1, & \mbox{if} & t=t_{0}\\
0, & \mbox{if} & t\in K\setminus B(t_{0},\delta).\end{array}\right.\]
 Since $x(t)>0$ only when $a_{0}(t)>0$, it follows that $x\in\overline{a_{0}Aa_{0}}$.
Moreover $x(t_{0})a_{0}(t_{0})x(t_{0})>1-\varepsilon$, hence $\left\Vert xa_{0}x\right\Vert >1-\varepsilon$.

Fix $t\in K$. If $t\in B(t_{0},\delta)$ then $\phi_{e}(t)\notin B(t_{0},\delta)$,
by our choice of $\delta_{1}$ and the fact that each map $\phi_{e}$
is a contraction, for all $e\in E$ such that $r(e)=v_{0}$; so $x\circ\phi_{e}(t)x(t)=0$.
If $t\notin B(t_{0},\delta)$, then $x(t)=0$, hence $x\circ\phi_{e}(t)x(t)=0$,
for all $e\in E$ such that $t\in K_{r(e)}$. Therefore $\left(\Phi(x)\xi x\right)(e,t)=x\circ\phi_{e}(t)\xi(e,t)x(t)=0$
for all $(e,t)\in E\times_{G}K$. Since\[
\langle\Phi(x)\xi x,\Phi(x)\xi x\rangle_{A}(t)=\sum_{\twolines{e\in E}{t\in K_{r(e)}}}(x\circ\phi_{e}(t))^{2}\mid\xi(e,t)\mid^{2}x(t)^{2},\]
 we see that $\left\Vert \Phi(x)\xi x\right\Vert =0$.

For $n_{0}=2$, we choose $t_{0}\in K$ such that $a_{0}(t_{0})>1-\varepsilon$
and $t_{0}$ is not a fixed point for any $\phi_{\alpha}$ with $\alpha\in E^{2}$.
Let $v_{0}\in V$ be such that $t_{0}\in K_{v_{0}}$. Let $\delta_{1}>0$
be such that $B(\phi_{\alpha}(t_{0}),\delta_{1})\cap B(t_{0},\delta_{1})=\emptyset$,
for all $\alpha\in E^{2}$ for which $r(\alpha)=v_{0}$, and such
that $B(t_{0},\delta_{1})\subset K_{v_{0}}$. Choosing $\delta_{2},\delta$
and $x$ as before, we conclude that $x\in\overline{a_{0}Aa_{0}}$
and $\left\Vert x\right\Vert >1-\varepsilon$. Moreover, we have $x\circ\phi_{\alpha}(t)x(t)=0$
for all $t\in K$, $\alpha\in E\cup E^{2}$ (since $\phi_{\alpha}$
is a contraction, for all $\alpha\in E\cup E^{2}$); and since\[
\left\langle \Phi^{(2)}(x)\xi^{2}x,\Phi^{(2)}(x)\xi^{2}x\right\rangle _{A}(t)=\sum_{\twolines{\alpha\in E^{2}}{t\in K_{r(\alpha)}}}(x\circ\phi_{\alpha}(t))^{2}\left|\xi_{2}^{2}(\alpha_{2},t)\right|^{2}\left|\xi_{1}^{2}(\alpha_{1},\phi_{\alpha_{1}}(t))\right|^{2}x(t)^{2}=0,\]
 it follows that $\left\Vert \Phi^{(k)}(x)\xi^{k}x\right\Vert =0$
for $k=1,2$. Applying the same argument inductively, we see that
$\mathcal{X}$ is an aperiodic $C^{*}$-correspondence. 
\end{proof}
Let $K^{1}$ and $K^{2}$ be two compact metric spaces. Let $A_{1}=C(K^{1})$
and $A_{2}=C(K^{2})$. If $\mequiv{A_{1}}{A_{2}}{\mathcal{Z}}$, then
the Rieffel correspondence determines a unique homeomorphism $f:K^{1}\rightarrow K^{2}$
and a unique Hermitian line bundle $\mathcal{L}$ over $\Graph(f)=\{(x,f(x))\;:\; x\in K^{1}\}$,
such that $\mathcal{Z}$ is isomorphic to $\Gamma(\mathcal{L})$ (see
\cite{Rieffel2}, \cite[Section 3.3 and Example 4.55]{Raburn-Williams},
\cite[Appendix (A)]{Raeburn-Picardgroup}), where $\Gamma(\mathcal{L})$
is the imprimitivity bimodule of the cross sections of $\mathcal{L}$
endowed with the following structure:\begin{eqnarray*}
(a\cdot s\cdot b)(x,f(x)) & = & a(x)s(x,f(x))b(f(x),\\
\langle s_{1},s_{2}\rangle_{A_{2}}(y) & = & \overline{s_{1}(f^{-1}(y),y)}s_{2}(f^{-1}(y),y),\\
\mbox{and}\;_{A_{1}}\langle s_{1},s_{2}\rangle(x) & = & s_{1}(x,f(x))\overline{s_{2}(x,f(x))},\end{eqnarray*}
for all $a\in A_{1}$, $b\in A_{2}$, $s,s_{1},s_{2}\in\Gamma(\mathcal{L})$.
We write $\mathcal{Z}(f,\mathcal{L})$ for $\Gamma(\mathcal{L})$.

We are ready to proof the main theorem.

\noindent \emph{Proof} (\emph{of Theorem} \ref{thm:Main}). By Proposition
\ref{pro:aperiodicity}, $\mathcal{X}_{1}$ and $\mathcal{X}_{2}$
are aperiodic $C^{\ast}$-correspondences. Using \cite[Theorem 7.2]{Muhly-Baruc-Moritaequiv},
we obtain that (\ref{enu:first}) implies (\ref{enu:second}).

Now we show that (\ref{enu:second}) implies (\ref{enu:third}). Suppose
that $\mathcal{X}_{1}$ and $\mathcal{X}_{2}$ are strongly Morita
equivalent. This implies that $A_{1}$ and $A_{2}$ are strongly Morita
equivalent via an imprimitivity bimodule $\mathcal{Z}$ such that
$\mathcal{Z}\otimes\mathcal{X}_{2}$ is isomorphic to $\mathcal{X}_{1}\otimes\mathcal{Z}$.
Let $f:K^{1}\rightarrow K^{2}$ and $\mathcal{L}$ be the homeomorphism
and the line bundle determined by the Rieffel correspondence. We have
that $\mathcal{Z}(f,\mathcal{L})\otimes\mathcal{X}_{2}$ is isomorphic
to $\mathcal{X}_{1}\otimes\mathcal{Z}(f,\mathcal{L})$. Hence $\mathcal{Z}(f,\mathcal{L})\otimes\mathcal{X}_{2}\otimes\widetilde{\mathcal{Z}(f,\mathcal{L})}$
is isomorphic to $\mathcal{X}_{1}$, where $\widetilde{\mathcal{Z}(f,\mathcal{L})}$
is the dual imprimitivity bimodule (see \cite[Proposition 3.18]{Raburn-Williams}).
We prove that $\mathcal{Z}(f,\mathcal{L})\otimes\mathcal{X}_{2}\otimes\widetilde{\mathcal{Z}(f,\mathcal{L})}$
is isomorphic to $\mathcal{X}_{2}$ over an isomorphism $\alpha$
of $A_{1}$ and $A_{2}$.

Let $\alpha:A_{1}\rightarrow A_{2}$ be defined by the formula $\alpha(a)=a\circ f^{-1}$
and let $V:\mathcal{Z}(f,\mathcal{L})\otimes\mathcal{X}_{2}\otimes\widetilde{\mathcal{Z}(f,\mathcal{L})}\rightarrow\mathcal{X}_{2}$
be defined by the formula\[
V(s_{1}\otimes\xi\otimes\widetilde{s_{2}})(e,y)=s_{1}(f^{-1}(\phi_{e}^{2}(y)),\phi_{e}^{2}(y))\xi(e,x)\overline{s_{2}(f^{-1}(y),y)}.\]
Then $\alpha$ is an isomorphism and\[
V(a\cdot s_{1}\otimes\xi\otimes\widetilde{s_{2}}\cdot b)=a\cdot V(s_{1}\otimes\xi\otimes\widetilde{s_{2}})\cdot b,\]
for all $a,b\in A$, $s_{1},s_{2}\in\mathcal{Z}(f,\mathcal{L})$,
$\xi\in\mathcal{X}_{2}$. Moreover we have that\begin{eqnarray*}
 &  & \langle V(s_{1}\otimes\xi\otimes\widetilde{s_{2}}),V(t_{1}\otimes\eta\otimes\widetilde{t_{2}}\rangle_{A_{2}}(y)\\
 &  & =\sum_{\twolines{e\in E}{y\in K_{r(e)}^{2}}}\overline{V(s_{1}\otimes\xi\otimes\widetilde{s_{2}})}(e,y)V(t_{1}\otimes\eta\otimes\widetilde{t_{2}})(e,y)\\
 &  & =\sum_{\twolines{e\in E}{y\in K_{r(e)}^{2}}}\left(\overline{s_{1}(f^{-1}(\phi_{e}^{2}(y)),\phi_{e}^{2}(y))\xi(e,x)\overline{s_{2}(f^{-1}(y),y)}}\right.\\
 &  & \left.\cdot t_{1}(f^{-1}(\phi_{e}^{2}(y)),\phi_{e}^{2}(y))\eta(e,x)\overline{t_{2}(f^{-1}(y),y)}\right)\\
 &  & =\langle s_{1}\otimes\xi\otimes s_{2},t_{1}\otimes\eta\otimes t_{2}\rangle_{A_{2}},\end{eqnarray*}
for all $s_{1},s_{2},t_{1},t_{2}\in\mathcal{Z}(f,\mathcal{L})$ and
$\xi,\eta\in\mathcal{X}_{2}$. Also, for $\xi\in\mathcal{X}_{2}$,
$V(1\otimes\xi\otimes1)=\xi$. Hence $V$ is a correspondence isomorphism.
Thus $\mathcal{X}_{1}$ is isomorphic to $\mathcal{X}_{2}$.

The rest is clear. \hfill \qed

It was shown in \cite[Theorem 2.3]{John} that the Cuntz-Pimsner algebra
of the $C^{*}$-correspondence built from a Mauldin-Williams graph
is isomorphic to the Cuntz-Krieger algebra of the underlying graph
$G=(V,E,r,s)$ (as defined in \cite{Kumjian-Pask-Raeburn-directedgraphs}).
Hence, for $C^{*}$-correspondences associated to Mauldin-Williams
graphs with the same underlying graph which are not isomorphic, we
obtain tensor algebras which are not Morita equivalent, but have the
same $C^{*}$-envelope, namely the Cuntz-Krieger algebra of the graph
$G$.

\section{The Isomorphism class of the $C^{*}$-Correspondences Associated
With Mauldin Williams Graphs}

In the following we analyze the relation between the isomorphism class
of the $C^{\ast}$-correspondences associated with two Mauldin-Williams
graphs, $\MWGs$, $i=1,2$ and the topological and dynamical properties
of the Mauldin-Williams graphs.

Since, by \cite[Section 4.2]{Pinzari-Watatami-Y} and \cite[Theorem 2.3]{John},
the Cuntz-Pimsner algebra associated to a Mauldin-Williams graph depends
only on the structure of the underlying graph $G$, we will consider
only Mauldin-Williams graphs having the same underlying graph $G=(V,E,r,s)$.

Next we determine necessary and sufficient conditions for the isomorphism
of the $C^{\ast}$-correspondences associated to two Mauldin-Williams
graphs.

\begin{prop}
\label{pro:SufficientME}For $i=1,2$, let $\mathcal{G}_{i}=(G,(K_{v}^{i})_{v\in V},(\phi_{e}^{i})_{e\in E})$
be two Mauldin-Williams graphs over the same underlying graph $G$.
Let $A_{i}=C(K^{i})$, $i=1,2$, be the associated $C^{*}$-algebras
and let $\mathcal{X}_{i}$, $i=1,2$, be the associated $C^{*}$-correspondences.
If there is a homeomorphism $f:K^{1}\rightarrow K^{2}$, a partition
of open subsets $\{ U_{1},\dots,U_{m}\}$ for $K^{1}$, for some $m\in\mathbb{N}$,
and if for each $U_{j}$ there is a permutation $\sigma_{j}\in S_{n}$,
where $n=|E|$, such that $f^{-1}\circ\phi_{\sigma_{j}(e)}^{2}\circ f|_{U_{j}}=\phi_{e}^{1}|_{U_{j}}$
and $f(K_{r(e)}^{1})=K_{r(\sigma_{j}(e))}^{2}$ for all $e\in E$,
$j\in\{1,\dots,m\}$, then $\mathcal{X}_{1}$ and $\mathcal{X}_{2}$
are isomorphic. 
\end{prop}
\begin{proof}
Since $f$ is a homeomorphism, the map $\beta:A_{2}\rightarrow A_{1}$,
defined by the equation $\beta(b)=b\circ f$ for all $b\in A_{2}$,
is a $C^{\ast}$-isomorphism. Define $V:\mathcal{X}_{2}\rightarrow\mathcal{X}_{1}$
by the formula\[
V(\xi)(e,x)=\sum_{k=1}^{m}\xi_{\sigma_{k}(e)}(f(x))\cdot1_{U_{k}}(x),\]
for all $(e,x)\in E\times_{G}K$, where $\xi_{\sigma_{k}(e)}(f(x)):=\xi(\sigma_{k}(e),f(x))$.
We show that $V$ is a $C^{\ast}$-correspondence isomorphism over
$\beta$. Let $b_{1},b_{2}\in A_{2}$ and $\xi\in\mathcal{X}_{2}$.
We have\begin{eqnarray*}
V(b_{1}\cdot\xi\cdot b_{2})(e,x) & = & \sum_{k=1}^{m}b_{1}\circ\phi_{\sigma_{k}(e)}^{2}(f(x))\xi_{\sigma_{k}(e)}(f(x))b_{2}(f(x))1_{U_{k}}(x)\\
 & = & \sum_{k=1}^{m}b_{1}\circ f\circ\phi_{e}^{1}(x)\xi_{\sigma_{k}(e)}(f(x))1_{U_{k}}(x)\cdot\beta(b_{2})(x)\\
 & = & \beta(b_{1})\cdot V(\xi)\cdot\beta(b_{2})(e,x).\end{eqnarray*}
Also\[
\langle V(\xi),V(\eta)\rangle_{A_{1}}(x)=\sum_{\genfrac{}{}{0pt}{}{e\in E}{f(x)\in K_{r(e)}^{2}}}\left(\sum_{k=1}^{m}\overline{\xi_{\sigma_{k}(e)}(f(x)}\eta_{\sigma_{k}(e)}(f(x))1_{U_{k}}(x)\right),\]
hence $\langle V(\xi),V(\eta)\rangle_{A_{1}}=\beta(\langle\xi,\eta\rangle_{A_{2}})$.
Finally one can see that $V$ is onto, hence $V$ is a $C^{\ast}$-correspondence
isomorphism.
\end{proof}
Recall that, for $k\geq2$, $E^{k}:=\{\alpha=(\alpha_{1},\cdots,\alpha_{k})\;:\;\alpha_{i}\in E\;\mbox{and}\; r(\alpha_{i})=s(\alpha_{i+1}),i=1,\cdots,k-1\}$,
is the set of \emph{paths of length $k$} in the graph $G$. Let $E^{\ast}=\bigcup_{k\in N}E^{k}$
be the space of finite paths in the graph $G$. Also the \emph{infinite
path space, $E^{\infty}$,} is defined to be \[
E^{\infty}:=\{(\alpha_{i})_{i\in\mathbb{N}}\quad:\quad\alpha_{i}\in E\quad\mbox{and}\quad r(\alpha_{i})=s(\alpha_{i+1})\;\mbox{for all}\: i\in\mathbb{N}\}.\]

\newcommand{\MWG}{\mathcal{M}=(G,\{ K_{v},\rho_{v}\}_{v\in E^{0}},\{\phi_{e}\}_{e\in E^{1}})}
For $v\in V$, we also define $E^{k}(v):=\{\alpha\in E^{k}\;:\; s(\alpha)=v\}$,
and $E^{\ast}(v)$ and $E^{\infty}(v)$ are defined similarly. We
consider $E^{\infty}(v)$ to be endowed with the metric: $\delta_{v}(\alpha,\beta)=c^{|\alpha\wedge\beta|}$
if $\alpha\ne\beta$ and $0$ otherwise, where $\alpha\wedge\beta$
is the longest common prefix of $\alpha$ and $\beta$, and $|w|$
is the length of the word $w\in E^{\ast}$ (see \cite[Page 116]{Edgar}).
Then $E^{\infty}(v)$ is a compact metric space, and, since $E^{\infty}$
equals the disjoint union of the spaces $E^{\infty}(v)$, $E^{\infty}$
becomes a compact metric space in a natural way. Define the maps $\theta_{e}:E^{\infty}(r(e))\rightarrow E^{\infty}(s(e))$
by the formula $\theta_{e}(\alpha)=e\alpha,\;\mbox{for all}\:\alpha\in E^{\infty}\mbox{ and for all}\: e\in E$.
Then $(G,(E^{\infty}(v))_{v\in V},(\theta_{e})_{e\in E})$ is a Mauldin-Williams
graph. We set $A_{E}:=C(E^{\infty})$ and we set $\mathcal{E}$ be
the $C^{*}$-correspondence associated to this Mauldin-Williams graph.
Let $\MWG$ be a Mauldin-Williams graph. For $(\alpha_{1},\dots,\alpha_{n})\in E^{n}$
let $K_{(\alpha_{1},\dots,\alpha_{n})}:=\phi_{\alpha_{1}}\circ\dots\phi_{\alpha_{n}}(K_{r(\alpha_{n})})$.
Then, for any infinite path $\alpha=(\alpha_{n})_{n\in\mathbb{N}}\in E^{\infty}$,
$\bigcap_{n\ge1}K_{(\alpha_{1},\dots,\alpha_{n})}$ contains only
one point. Therefore we can define a map $\pi:E^{\infty}\rightarrow K$
by $\{\pi(x)\}=\bigcap_{n\ge1}K_{(\alpha_{1},\dots,\alpha_{n})}$.
Since $\pi(E^{\infty})$ is also an invariant set, $\pi$ is a continuous,
onto map and $\pi(E^{\infty}(v))=K_{v}$. Moreover, $\pi\circ\theta_{e}=\phi_{e}\circ\pi$. 

We say that a Mauldin-Williams graph $\MWG$ is \emph{totally disconnected}
if $\phi_{e}(K_{r(e)})\bigcap\phi_{f}(K_{r(f)})=\emptyset$ if $s(e)=s(f)$
and $e\ne f$.

\begin{cor}
Let $\MWG$ be a totally disconnected Mauldin-Williams graph. Let
$A$ be the $C^{*}$-algebra and $\mathcal{X}$ be the $C^{*}$-correspondence
associated to this Mauldin-Williams graph. Then $\mathcal{X}$ is
isomorphic with $\mathcal{E}$, as $C^{\ast}$-correspondences. In
particular, one obtains that for any two totally disconnected Mauldin-Williams
graphs having the same underlying graph $G$, the $C^{\ast}$-correspondences
and tensor algebras associated are isomorphic.
\end{cor}
\begin{proof}
If the Mauldin-Williams graph is totally disconnected, then the map
$\pi:E^{\infty}\rightarrow K$ defined above is a homeomorphism. Moreover
$\pi\circ\theta_{e}\circ\pi^{-1}=\phi_{e}$ for all $e\in E$, therefore
the associated $C^{\ast}$-correspondences are isomorphic. 
\end{proof}
The converse of this corollary is true and will be proved later.

The next theorem is a converse of the Proposition \ref{pro:SufficientME}.
We note, however, that the family of open sets $\{ U_{i}\}$ here
is not required to be a partition of the compact set $K^{1}$, but
only a finite open cover of it.

\begin{thm}
\label{thmequivofMWg}For $i=1,2$, let $\mathcal{G}_{i}=(G,(K_{v}^{i})_{v\in V},(\phi_{e}^{i})_{e\in E})$
be two Mauldin-Williams graphs over the same underlying graph $G$.
Let $A_{i}=C(K^{i})$, $i=1,2$, be the associated $C^{*}$-algebras
and let $\mathcal{X}_{i}$, $i=1,2$, be the associated $C^{*}$-correspondences.
If $\mathcal{X}_{1}$ and $\mathcal{X}_{2}$ are isomorphic, then
there is a homeomorphism $f:K^{1}\rightarrow K^{2}$, a finite open
cover of $K^{1}$, $\{ U_{1},\dots,U_{m}\}$, and for each $U_{j}$
there is a permutation $\sigma_{j}\in S_{n}$ ($n=|E|$) such that
\begin{equation}
f^{-1}\circ\phi_{e}^{2}\circ f|_{U_{j}}=\phi_{\sigma_{j}(e)}^{1}|_{U_{j}}\;\mbox{for all}\: e\in E,i\in\{1,\dots,m\}.\label{eq:equivofMWg}\end{equation}

\end{thm}
\begin{proof}
Since $\mathcal{X}_{1}$ and $\mathcal{X}_{2}$ are isomorphic, there
is a $C^{\ast}$-isomorphism $\beta:A_{2}\rightarrow A_{1}$ and a
$C^{\ast}$-correspondence isomorphism $W:\mathcal{X}_{2}\rightarrow\mathcal{X}_{1}$
such that $W(b_{1}\cdot\xi\cdot b_{2})=\beta(b_{1})W(\xi)\beta(b_{2})$
and $\langle W(\xi),W(\eta)\rangle_{A_{1}}=\beta(\langle\xi,\eta\rangle_{A_{2}})$,
for all $b_{1},b_{2}\in A_{2}$, $\xi,\eta\in\mathcal{X}_{2}$. Let
$f:K^{1}\rightarrow K^{2}$ be the homeomorphism which implements
$\beta$, that is $\beta(b)=b\circ f$, for all $b\in A_{2}$.

Let $\delta_{e}\in\mathcal{X}_{2}$, defined by\[
\delta_{e}(g,y)=\left\{ \begin{array}{rcl}
1, & \mbox{if }e=g\\
0, & \mbox{otherwise},\end{array}\right.\]
for $e\in E$, be the natural basis in $\mathcal{X}_{2}$ and let
$\left(\delta_{e}^{\prime}\right)_{e\in E}\subset\mathcal{X}_{1}$
be the natural basis in $\mathcal{X}_{1}$, which is defined similarly.

For $\xi\in\mathcal{X}_{2}$, $\xi=\sum_{g\in E}\delta_{g}\cdot\xi_{g}$,
where $\xi_{g}(y)=\xi(g,y)$ for all $y\in K_{r(g)}^{2}$ and is $0$
otherwise. With respect to the bases, we can write\begin{align}
W(\xi)= & W\left(\sum_{g\in E}\delta_{g}\cdot\xi_{g}\right)=\sum_{g\in E}W(\delta_{g})\cdot\xi_{g}\circ f\nonumber \\
= & \sum_{g\in E}\sum_{e\in E}\delta_{e}^{\prime}\cdot w_{eg}\xi_{g}\circ f,\label{eq:isomorform}\end{align}
 where\begin{equation}
W(\delta_{g})=\sum_{e\in E}\delta_{e}^{\prime}\cdot w_{eg}\quad,w_{eg}\in A_{1}\label{eq:isomform2}\end{equation}
 and $w_{eg}$ are given by the formula $w_{eg}=\langle\delta_{e}^{\prime},W(\delta_{g})\rangle_{A_{1}}$,
for all $e,g\in E$. We call $(w_{eg})_{e,g\in E}$ the matrix of
$W$ with respect to the basis $(\delta_{e}^{\prime})_{e\in E}$ and
$(\delta_{g})_{g\in E}$ (it is an $n\times n$ matrix, where $n=|E|$).
Since $W$ preserves the inner product, we see that\begin{equation}
\langle W(\delta_{g}),W(\delta_{e})\rangle=\langle\delta_{g},\delta_{e}\rangle=\delta_{ge},\label{eq:iso1}\end{equation}
 where $\delta_{ge}(x)=1$ if $e=g$ and $x\in K_{r(e)}$ and is $0$
otherwise. Also\begin{align}
\langle W(\delta_{g}),W(\delta_{e})\rangle= & \sum_{f\in E}w_{fg}^{*}w_{fe}.\label{iso2}\end{align}
 Equations (\ref{eq:iso1}) and (\ref{iso2}) imply that for every
$x\in K^{1}$ the matrix $(w_{ef}(x))_{e,f\in E}$ is invertible.
Hence there is $\sigma_{x}\in S_{n}$ such that $w_{\sigma_{x}(e)e}(x)\neq0$
for all $e\in E$. Therefore there is a neighborhood $U_{x}$ of $x$
such that\begin{equation}
w_{\sigma_{x}(e)e}(y)\neq0\;\mbox{for all}\,\, e\in E,y\in U_{x}\;\mbox{and}\; x\in K^{1}.\label{eq:neigbhoodforx}\end{equation}
Let $b\in A_{2}$. Then, for $h\in E$ we have that\[
W(b\cdot\delta_{h})=\sum_{e\in E}\delta_{e}^{\prime}w_{eh}b\circ\phi_{h}^{2}\circ f\]
and\[
\beta(b)\cdot W(\delta_{h})=\sum_{e\in E}\delta_{e}^{\prime}b\circ f\circ\phi_{e}^{1}w_{eh}.\]
 Fix $x\in K^{1}$ and let $\sigma_{x}\in S_{n}$ and $U_{x}$ be
defined as in Equation (\ref{eq:neigbhoodforx}). Then\[
W(b\cdot\delta_{h})(\sigma_{x}(h),y)=w_{\sigma_{x}(h)h}(y)b\circ\phi_{h}^{2}\circ f(y)\]
and\[
(\beta(b)\cdot W(\delta_{h}))(\sigma_{x}(h),y)=b\circ f\circ\phi_{\sigma_{x}(h)}^{1}(y)w_{\sigma_{x}(h)h}(y)\]
 for all $y\in U_{x}$ and for all $h\in E$. Since $W$ is a $C^{\ast}$-correspondence
isomorphism and $w_{\sigma_{x}(h)h}(y)\neq0$ for all $y\in U_{x}$,
for any $x\in K^{1}$, there is a neighborhood $U_{x}$ of $x$ in
$K^{1}$ and there is a permutation $\sigma_{x}\in S_{n}$ such that\[
f^{-1}\circ\phi_{h}^{2}\circ f|_{U_{x}}=\phi_{\sigma_{x}(h)}^{1}|_{U_{x}}\;\mbox{for all}\: h\in E.\]
 Hence we can find a finite cover $\{ U_{1},\dots,U_{m}\}$ of $K^{1}$
and for each $U_{i}$ we can find a permutation $\sigma_{i}\in S_{n}$
such that the Equation (\ref{eq:equivofMWg}) holds.
\end{proof}
In the special case when the two Mauldin-Williams graphs are totally
disconnected, more can be said about the choice of the permutations
$\sigma_{i}$.

\begin{cor}
\label{cor:IFStotdisc}Let $\MWGs$ be two Mauldin-Williams graphs.
Let $A_{i}=C(K^{i})$ and let $\mathcal{X}_{i}$, $i=1,2$, be the
associated $C^{*}$-algebras and $C^{*}$-correspondences. If $\mathcal{G}_{1}$
is totally disconnected and if $\mathcal{X}_{1}$ is isomorphic with
$\mathcal{X}_{2}$ there is a continuous map $h:K^{1}\rightarrow S_{n}$
such that $f^{-1}\circ\phi_{e}^{2}\circ f(x)=\phi_{h(x)(e)}(x),\;\mbox{for all}\; x\in K^{1}$. 
\end{cor}
\begin{proof}
Recall that if $\mathcal{G}_{1}$ is totally disconnected, then $\phi_{e}^{1}(K_{r(e)})\bigcap\phi_{f}^{1}(K_{r(f)})=\emptyset$
if $e\ne f$. From the Theorem \ref{thmequivofMWg} we know that there
are open sets $\{ U_{1},\cdots,U_{m}\}$, for some $m\in\mathbb{N}$,
and permutations $\sigma_{1},\cdots,\sigma_{m}\in S_{n}$ such that
\begin{equation}
f^{-1}\circ\phi_{e}^{2}\circ f|_{U_{i}}=\phi_{\sigma_{i}(e)}^{1}|_{U_{i}}\quad\mbox{for all}\: e\in E,i\in\{1,\cdots,m\}.\label{eq:equivfortotdisconn1}\end{equation}
 If $U_{i}\cap U_{j}\ne\emptyset$ for some $i\ne j$, then it follows
that $\phi_{\sigma_{i}(e)|U_{i}\cap U_{j}}^{1}=\phi_{\sigma_{j}(e)}^{1}|_{U_{i}\cap U_{j}}$
for all $e\in E$, hence $\sigma_{i}(e)=\sigma_{j}(e)\;\mbox{for all}\: e\in E$,
so $\sigma_{i}=\sigma_{j}$. Therefore we can choose the cover $U_{1},\cdots,U_{m}$
such that $U_{i}\cap U_{j}=\emptyset$ if $i\ne j$.

Let $x\in K^{1}$. Then there is a unique $i\in\{1,\cdots,n\}$ such
that $x\in U_{i}$. We define $h(x)=\sigma_{i}$. Then $h:K^{1}\rightarrow S_{n}$
is a well defined map. Moreover, $h$ is continuous (considering $S_{n}$
endowed with the discrete topology), since for every $\sigma\in S_{n}$,
$h^{-1}(\sigma)=\emptyset$ or $h^{-1}(\sigma)=U_{i}$, for some $i\in\{1,\cdots,n\}$.
Finally, from the Equation (\ref{eq:equivfortotdisconn1}) we obtain
that \[
f^{-1}\circ\phi_{e}^{2}\circ f(x)=\phi_{h(x)(e)}^{1}(x)\;\mbox{for all}\: x\in K^{1}\;\mbox{and}\: e\in E.\]

\end{proof}
Suppose that $\MWGs$ are two Mauldin-Williams graphs, that satisfy
the hypothesis of the Corollary \ref{cor:IFStotdisc}. We claim that
$\mathcal{G}_{2}$ must also be totally disconnected. Suppose that
there are $e,f\in E$, $e\ne f$, such that $\phi_{e}^{2}(K_{r(e)}^{2})\bigcap\phi_{f}^{2}(K_{r(f)}^{2})\ne\emptyset$.
Then there is an $x\in K^{1}$ such that $y:=f(x)\in\phi_{e}^{2}(K_{r(e)}^{2})\bigcap\phi_{f}^{2}(K_{r(f)}^{2})$.
Then $\phi_{h(x)(e)}^{1}(x)=\phi_{h(x)(f)}^{1}(x)$, which is a contradiction
since $\mathcal{G}_{1}$ is totally disconnected. So $\mathcal{G}_{2}$
is totally disconnected.

\begin{example}
Let $K$ be the Cantor set, let $\phi_{i}:K\rightarrow K$, $i=1,2$
be the maps defined by the formulae $\phi_{1}(x)=\frac{1}{3}x$ and
$\phi_{2}(x)=\frac{1}{3}x+\frac{2}{3}$. Then $K$ is the invariant
set of $(\phi_{1},\phi_{2})$. Let $T=[0,1]$ and let $\psi_{i}:T\rightarrow T$,
$i=1,2$ be the maps defined by the formulae $\psi_{1}(x)=\frac{1}{2}x$
and $\psi_{2}(x)=-\frac{1}{2}x+1$. Then $T$ is the invariant set
of $(\psi_{1},\psi_{2})$. Since $(\psi_{1},\psi_{2})$ is not totally
disconnected, we see that the associated $C^{\ast}$-correspondences
are not strongly Morita equivalent. Hence the tensor algebras fail
to be strongly Morita equivalent in the sense of \cite{Blecher-Muhly-Paulsen},
yet their $C^{\ast}$-envelopes coincide with $\mathcal{O}_{2}$. 
\end{example}

\begin{example}
Let $T$ be the regular triangle in $\mathbb{R}^{2}$ with vertices
$A=(0,0)$, $B=(1,0)$ and $C=(1/2,\sqrt{3}/2)$. Let $\phi_{1}(x,y)=\left(\frac{x}{2}+\frac{1}{4},\frac{y}{2}+\frac{\sqrt{3}}{4}\right)$,
$\phi_{2}(x,y)=\left(\frac{x}{2},\frac{y}{2}\right)$ and $\phi_{3}(x,y)=\left(\frac{x}{2}+\frac{1}{2},\frac{y}{2}\right)$.
Then the invariant set $K$ of this iterated function system is the
Sierpinski gasket. Let $\psi_{1}=\sigma_{1}\circ\phi_{1}$, $\psi_{2}=\phi_{2}$
and $\psi_{3}=\sigma_{3}\circ\phi_{3}$, where $\sigma_{i}$ is the
symmetry about the median from the vertex $\phi_{i}(C)$ of the triangle
$\phi_{i}(T)$, for $i=1,3$. Then the invariant set of this iterated
function system is also the Sierpinski gasket, but the $C^{\ast}$-correspondences
associated to $(\phi_{1},\phi_{2},\phi_{3})$ and $(\psi_{1},\psi_{2},\psi_{3})$
are not isomorphic.
\end{example}

\end{document}